\documentclass[12pt]{amsart}

\usepackage{mathrsfs}

\newtheorem{theorem}{Theorem}[section]

\theoremstyle{definition}

\theoremstyle{remark}
\newtheorem{remark}[theorem]{Remark}

\numberwithin{equation}{section}

\let \la=\lambda
\let \e=\varepsilon
\let \d=\delta
\let \o=\omega
\let \a=\alpha
\let \f=\varphi
\let \b=\beta

\let \g=\gamma
\let \O=\Omega
\let \si=\sigma

\let \G=\Gamma
\let \ga=\gamma

\begin{document}
\title[Mixed $A_p$-$A_r$ inequalities]
{Mixed $A_p$-$A_r$ inequalities for classical singular integrals and
Littlewood-Paley operators}

\author{Andrei K. Lerner}
\address{Department of Mathematics,
Bar-Ilan University, 52900 Ramat Gan, Israel}
\email{aklerner@netvision.net.il}

\begin{abstract}
We prove mixed $A_p$-$A_r$ inequalities for several basic singular
integrals, Littlewood-Paley operators, and the vector-valued maximal
function. Our key point is that $r$ can be taken arbitrary big.
Hence such inequalities are close in spirit to those obtained
recently in the works by T. Hyt\"onen and C. P\'erez, and M. Lacey.
On one hand, the ``$A_p$-$A_{\infty}$" constant in these works
involves two independent suprema. On the other hand, the
``$A_p$-$A_r$" constant in our estimates involves a joint supremum
but of a bigger expression. We show in simple examples that both
such constants are incomparable. This leads to a natural conjecture
that the estimates of both types can be further improved.
\end{abstract}



\maketitle

\section{Introduction}
Given a weight (that is, a non-negative locally integrable function)
~$w$ and a cube $Q\subset {\mathbb R}^n$, let
$$A_p(w;Q)=\Big(\frac{1}{|Q|}\int_Qw\Big)\Big(\frac{1}{|Q|}
\int_Qw^{-\frac{1}{p-1}}\Big)^{p-1}\quad (1<p<\infty)$$ and
$$\|w\|_{A_p}=\sup_{Q\subset {\mathbb R}^n}A_p(w;Q).$$

Sharp weighted norm inequalities in terms of
$\|w\|_{A_p}$ have been obtained recently for the Calder\'on-Zygmund
operators and for a large class of the Littlewood-Paley operators.
To be more precise, if $T$ is a Calder\'on-Zygmund operator, then
\begin{equation}\label{cz}
\|T\|_{L^p(w)}\le
c(T,p,n)\|w\|_{A_p}^{\max(1,\frac{1}{p-1})}\quad(1<p<\infty).
\end{equation}
This result in its full generality is due to T. Hyt\"onen \cite{H};
we also refer to this work for a very detailed history of closely
related results and particular cases. Soon after appearing \cite{H},
a somewhat simplified approach to (\ref{cz}) was found in
\cite{HPTV}.

If $S$ is a Littlewood-Paley operator (in particular, any typical square function),
then (see \cite{L2} and the references therein)
\begin{equation}\label{lp}
\|S\|_{L^p(w)}\le
c(S,p,n)\|w\|_{A_p}^{\max(\frac{1}{2},\frac{1}{p-1})}\quad(1<p<\infty).
\end{equation}

Observe that the exponents in (\ref{cz}) and (\ref{lp}) are sharp
for any $1<p<\infty$. However, it turns out that this is not the end
of the story. Very recently, T. Hyt\"onen and C.
P\'erez \cite{HP} have studied mixed $A_p$-$A_{\infty}$ estimates that
improve many of known sharp $A_p$ estimates. Denote
$$\|w\|_{A_{\infty}}=\sup_{Q\subset {\mathbb R}^n}A_{\infty}(w;Q)=
\sup_{Q\subset {\mathbb R}^n}\Big(\frac{1}{|Q|}\int_Qw\Big)\exp\Big(\frac{1}{|Q|}\int_Q\log
w^{-1}\Big).$$ Set also
$$\|w\|_{A_{\infty}}'=\sup_{Q\subset {\mathbb R}^n}\frac{1}{w(Q)}\int_QM(w\chi_Q),$$
where $M$ is the Hardy-Littlewood maximal operator. Observe that
$$c_n\|w\|_{A_{\infty}}'\le \|w\|_{A_{\infty}}\le \|w\|_{A_p}\quad
(1<p<\infty)$$ and the first inequality here cannot be reversed (see
\cite{HP} for the details).

One of the main results in \cite{HP} is the following improvement of
(\ref{cz}) in the case $p=2$:
\begin{equation}\label{l2}
\|T\|_{L^2(w)}\le
c(T,n)\|w\|_{A_2}^{1/2}\max(\|w\|_{A_{\infty}}',\|w^{-1}\|_{A_{\infty}}').
\end{equation}
It is well known that the case $p=2$ is crucial for inequality
(\ref{cz}). Indeed, (\ref{cz}) for any $p\not=2$ follows from the
linear $L^2(w)$ bound and the sharp version of the Rubio de Francia
extrapolation theorem. Adapting such approach, the authors in
\cite{HP} extended (\ref{l2}) for any $p\not=2$. For example, it was
shown that for $p>2$,
\begin{equation}\label{lpw}
\|T\|_{L^p(w)}\le
c\|w\|_{A_p}^{\frac{2}{p}-\frac{1}{2(p-1)}}\big(\|w\|_{A_{\infty}}^{\frac{1}{2(p-1)}}+
\|\si\|_{A_{\infty}}^{\frac{1}{2}}\big)(\|w\|_{A_{\infty}}')^{1-\frac{2}{p}},
\end{equation}
where $\si=w^{-\frac{1}{p-1}}$.

It turns out that while the extrapolation method is powerful for
(\ref{cz}), it is not so effective for mixed $A_p$-$A_{\infty}$
inequalities. Indeed, T.~Hyt\"onen et al. \cite{HLMORSU} improved
(\ref{lpw}) (at least for $p>4$) without the use of extrapolation, namely, it is proved
in \cite{HLMORSU} that
\begin{equation}\label{imp}
\|T^*\|_{L^p(w)}\le
c(T,p,n)\Big(\|w\|_{A_p}^{1/p}(\|w\|_{A_{\infty}}')^{1/p'}+\|w\|_{A_p}^{\frac{1}{p-1}}\Big)\quad(p>1),
\end{equation}
where $T^*$ is the maximal Calder\'on-Zygmund operator.

Soon after that, M. Lacey \cite{L} improved (\ref{imp}) and (\ref{lpw}) for several
classical singular integrals:
\begin{equation}\label{Lac}
\|T^*\|_{L^p(w)}\le
c(T,p,n)\|w\|_{A_p}^{1/p}\max\Big(\|w\|_{A_{\infty}}')^{1/p'}
,(\|\si\|_{A_{\infty}}')^{1/p}\Big),
\end{equation}
and it was conjectured in \cite{L}
that (\ref{Lac}) holds for any Calder\'on-Zygmund operator.

More precisely, (\ref{Lac}) was proved for the Hilbert, Riesz and
Beurling operators and for any one-dimensional convolution
Calder\'on-Zygmund operator with odd $C^2$ kernel. All these
operators are unified by the fact that they can be represented as a
suitable average of the so-called Haar shift operators ${\mathbb S}$
with bounded complexity. In order to handle such operators, it was
used in \cite{L} a ``local mean oscillation" decomposition. The
latter decomposition was obtained by the author in \cite{L1}. Then,
its various applications (in particular, to the Haar shift
operators) have been found by D. Cruz-Uribe, J. Martell and C.
P\'erez in \cite{CMP}.

After an application of the decomposition to ${\mathbb S}$, the
proof of (\ref{Lac}) is reduced to showing that this estimate is
true for
$${\mathcal A}_{\gamma}f(x)=\sum_{j,k}\Big(\frac{1}{|\ga
Q_j^k|}\int_{\ga Q_j^k}|f|\Big)\chi_{Q_j^k}(x),$$ where $Q_j^k$ are
the dyadic cubes with good overlapping properties. This is done in
\cite{L} by means of a number of interesting tricks. It is mentioned
in \cite{L} that a more elementary approach to $A_{\ga}$ (used in
\cite{CMP} in order to prove (\ref{cz}) for classical singular
operators mentioned above) does not allow to get (\ref{Lac}).

In this paper we show, however, that a variation of the
approach to $A_{\ga}$ from \cite{CMP} allows to get mixed estimates
of a different type, namely, we obtain $L^p(w)$ bounds in terms of
$$\|w\|_{(A_p)^{\a}(A_r)^{\b}}=\sup_{Q\subset {\mathbb R}^n}A_p(w;Q)^{\a}A_r(w;Q)^{\b}$$
for suitable $\a$ and $\b$. The key point in our results below is
that $r$ can be taken arbitrary big (but with the implicit constant
growing exponentially in $r$). Therefore, our estimates can be also
considered as a kind of $A_p$-$A_{\infty}$ estimates. An important
feature of the expression defining $\|w\|_{(A_p)^{\a}(A_r)^{\b}}$ is
that only one supremum is involved. We will show in simple examples
that $\|w\|_{(A_p)^{\a}(A_r)^{\b}}$ is incomparable with the
right-hand side of (\ref{Lac}), that is, each of such expressions
can be arbitrary larger than the other. This fact indicates that the
estimates of both types can be further improved.

In the next theorem we suppose that $T^*$ is the same operator as in
(\ref{Lac}), namely,
$$T^*f(x)=\sup_{\e<\d}\Big|\int_{\e<|x-y|<\d}f(y)K(x-y)dy\Big|,$$
where $K$ is one of the following kernels: (i)
$K(x)=\frac{1}{x},n=1$; (ii) $k(x)=\frac{x_j}{|x|^{n+1}},n\ge 2$;
(iii) $K(z)=\frac{1}{z^2}, z\in {\mathbb C}$; (iv) $K(x)$ is any
odd, one-dimensional $C^2$ kernel satisfying $|K^{(i)}(x)|\le
c|x|^{-1-i}\quad(i=0,1,2).$

\begin{theorem}\label{mixcz} For any $2\le p\le r<\infty$,
$$\|T^*\|_{L^p(w)}\le
c(T,p,r,n)\|w\|_{(A_p)^{\frac{1}{p-1}}(A_r)^{1-\frac{1}{p-1}}}.$$
\end{theorem}

A similar result holds for the Littlewood-Paley operators satisfying
(\ref{cz}). In the next theorem, $S$ is either the dyadic square
function or the intrinsic square function (and hence the theorem is
also true for the Lusin area integral~$S(f)$, the Littlewood-Paley
function $g(f)$, the continuous square functions $S_{\psi}(f)$ and
$g_{\psi}(f)$).

\begin{theorem}\label{mixlp} For any $3\le p\le r<\infty$,
$$\|S\|_{L^p(w)}\le
c(S,p,r,n)\|w\|_{(A_p)^{\frac{1}{p-1}}(A_r)^{\frac{1}{2}-\frac{1}{p-1}}}.$$
\end{theorem}

Also, the result of the same type holds for the vector-valued
maximal function (see Remark \ref{vvm} below).

In Section 4, we show the sharpness of the exponent $\frac{1}{p-1}$
in Theorems~\ref{mixcz} and \ref{mixlp}. Also we show that the
right-hand side in Theorem~\ref{mixcz} is incomparable with the one
in (\ref{Lac}).

A natural question appearing here is whether the right-hand side in
Theorem~\ref{mixcz} can be replaced by
$$\|w\|_{(A_p)^{\frac{1}{p-1}}(A_{\infty})^{1-\frac{1}{p-1}}}=
\sup_{Q}A_p(w;Q)^{\frac{1}{p-1}}A_{\infty}(w;Q)^{1-\frac{1}{p-1}}$$
or by
$$\|w\|_{(A_p)^{\frac{1}{p-1}}(A'_{\infty})^{1-\frac{1}{p-1}}}=
\sup_{Q}A_p(w;Q)^{\frac{1}{p-1}}A'_{\infty}(w;Q)^{1-\frac{1}{p-1}},$$
where $A'_{\infty}(w;Q)=\frac{1}{w(Q)}\int_QM(w\chi_Q)$.

\section{Preliminaries}
\subsection{Haar shift operators}
Given a general dyadic grid ${\mathscr{D}}$ and $m,k\in {\mathbb
N}$, we say that ${\mathbb S}$ is a (generalized) Haar shift
operator with parameters $m,k$ if
$${\mathbb S}f(x)={\mathbb S}_{{\mathscr{D}}}^{m,k}f(x)=
\sum_{Q\in {\mathscr{D}}}\sum_{{Q',Q''\in
{\mathscr{D}},Q',Q''\subset Q}\atop
{\ell(Q')=2^{-m}\ell(Q),\ell(Q'')=2^{-k}\ell(Q)}}\frac {\langle
f,h_{Q'}^{Q''}\rangle}{|Q|}h_{Q''}^{Q'}(x),
$$
where $\ell(Q)$ is the side length of $Q$, $h_{Q'}^{Q''}$ is a (generalized) Haar function on $Q'$, and
$h_{Q''}^{Q'}$ is one on $Q''$ such that
$$\|h_{Q'}^{Q''}\|_{L^{\infty}}\|h_{Q''}^{Q'}\|_{L^{\infty}}\le 1.$$
The number $\max(m,k)$ is called the complexity of ${\mathbb S}$.

We refer to \cite{H} for a more detailed explanation of this
definition. Also, it is shown in \cite{H} that any
Calder\'on-Zygmund operator can be represented as a suitable average
of ${\mathbb S}_{{\mathscr{D}}}^{m,k}$ with respect to all diadic
grids ${\mathscr{D}}$ and all $m,k \in {\mathbb N}$. In the case of
the classical convolution operators mentioned in Theorem
\ref{mixcz}, such an average can be taken only of  ${\mathbb
S}_{{\mathscr{D}}}^{m,k}$ with bounded complexity. This fact was
proved in the works \cite{DV} (the Beurling operator), \cite{P} (the
Hilbert transform), \cite{PTV} (the Riesz transforms), \cite{V} (any
one-dimensional singular integral with odd $C^2$ kernel).

Similarly to the maximal singular integral $T^*$, one can define the
maximal Haar shift operator ${\mathbb S}^{*}$, and to get a control
of $T^*$ by ${\mathbb S}^{*}$ (see \cite[Prop. 2.8]{HLMORSU}). In
particular, it suffices to prove Theorem \ref{mixcz} for a single
${\mathbb S}^*$ instead of $T^*$.

\subsection{Littlewood-Paley operators}
The dyadic square function is defined by
$$S_df(x)=\left(\sum_{Q\in{\mathcal D}}(f_Q-f_{\widehat
Q})^2\chi_Q(x)\right)^{1/2},$$ where the sum is taken over all
dyadic cubes on ${\mathbb R}^n$.

Let ${\mathbb R}^{n+1}_+={\mathbb R}^{n}\times{\mathbb R}_{+}$ and
$\G(x)=\{(y,t)\in {\mathbb{R}}^{n+1}_+:|y-x|<t\}$. For
$0<\a\le 1$, let ${\mathcal C}_{\a}$ be the family of functions
supported in $\{x:|x|\le 1\}$, satisfying $\int\psi=0$, and such
that for all $x$ and $x'$, $|\f(x)-\f(x')|\le |x-x'|^{\a}$. If $f\in
L^1_{\text{loc}}({\mathbb R}^n)$ and $(y,t)\in {\mathbb R}^{n+1}_+$,
we define
$$
A_{\a}(f)(y,t)=\sup_{\f\in {\mathcal C}_{\a}}|f*\f_t(y)|.
$$

The intrinsic square function $G_{\a}(f)$ is defined by
$$
G_{\a}(f)(x)=\left(\int_{\G(x)}
\big(A_{\a}(f)(y,t)\big)^2\frac{dydt}{t^{n+1}}\right)^{1/2}.
$$
This operator was introduced by M. Wilson \cite{W}.
On one hand $G_{\a}$ pointwise
dominates the classical and continuous $S$ and $g$ functions. On the
other hand, it is not essentially larger than any one of them.

Denote
$$T(Q)=\{(y,t)\in {\mathbb R}^{n+1}_+:y\in Q, \ell(Q)/2\le t<\ell(Q)\}$$
and $\g_Q(f)^2=\int_{T(Q)}
\big(A_{\a}(f)(y,t)\big)^2\frac{dydt}{t^{n+1}}$, and let
$$
\widetilde G_{\a}(f)(x)=\Big(\sum_{Q\in {\mathcal
D}}\g_Q(f)^2\chi_{3Q}(x)\Big)^{1/2}.
$$
Then we have that (see \cite{L2})
$$G_{\a}(f)(x)\le \widetilde G_{\a}(f)(x)\le c(\a,n)G_{\a}(f)(x).$$

\subsection{A ``local mean oscillation" decomposition}
Given a measurable function $f$ on ${\mathbb R}^n$ and a cube $Q$,
define the local mean oscillation of $f$ on $Q$ by
$$\o_{\la}(f;Q)=\inf_{c\in {\mathbb R}}
\big((f-c)\chi_{Q}\big)^*\big(\la|Q|\big)\quad(0<\la<1),$$ where
$f^*$ denotes the non-increasing rearrangement of $f$.

By a median value of $f$ over $Q$ we mean a possibly nonunique, real
number $m_f(Q)$ such that
$$\max\big(|\{x\in Q: f(x)>m_f(Q)\}|,|\{x\in Q: f(x)<m_f(Q)\}|\big)\le |Q|/2.$$

Given a cube $Q_0$, denote by ${\mathcal D}(Q_0)$ the set of all
dyadic cubes with respect to $Q_0$. If $Q\in {\mathcal D}(Q_0)$ and
$Q\not=Q_0$, we denote by $\widehat Q$ its dyadic parent, that is,
the unique cube from ${\mathcal D}(Q_0)$ containing $Q$ and such
that $|\widehat Q|=2^n|Q|$.

The dyadic local sharp maximal
function $M^{\#,d}_{\la;Q_0}f$ is defined by
$$M^{\#,d}_{\la;Q_0}f(x)=\sup_{x\in Q'\in
{\mathcal D}(Q_0)}\o_{\la}(f;Q').$$

The following theorem was proved in \cite{L1}.

\begin{theorem}\label{decom} Let $f$ be a measurable function on
${\mathbb R}^n$ and let $Q_0$ be a fixed cube. Then there exists a
(possibly empty) collection of cubes $Q_j^k\in {\mathcal D}(Q_0)$
such that
\begin{enumerate}
\renewcommand{\labelenumi}{(\roman{enumi})}
\item
for a.e. $x\in Q_0$,
$$
|f(x)-m_f(Q_0)|\le
4M_{1/4;Q_0}^{\#,d}f(x)+4\sum_{k=1}^{\infty}\sum_j
\o_{\frac{1}{2^{n+2}}}(f;\widehat Q_j^k)\chi_{Q_j^k}(x);
$$
\item
for each fixed $k$ the cubes $Q_j^k$ are pairwise disjoint;
\item
if $\O_k=\cup_jQ_j^k$, then $\O_{k+1}\subset \O_k$;
\item
$|\O_{k+1}\cap Q_j^k|\le \frac{1}{2}|Q_j^k|.$
\end{enumerate}
\end{theorem}

We shall use below the standard fact following from the above
properties (ii)-(iv), namely, that the sets $E_j^k=Q_j^k\setminus
\O_{k+1}$ are pairwise disjoint and $|E_j^k|\ge \frac{1}{2}|Q_j^k|$.

\section{Proof of Theorems \ref{mixcz} and \ref{mixlp}}
The key result implying both Theorems \ref{mixcz} and \ref{mixlp}
can be described as follows.

\begin{theorem}\label{key} Let $T$ be a sublinear operator
satisfying
\begin{equation}\label{cond}
\o_{\la}(|Tf|^{\nu};Q)\le c\Big(\frac{1}{|\ga Q|}\int_{\ga
Q}|f|dx\Big)^{\nu}
\end{equation}
for any dyadic cube $Q\subset {\mathbb R}^n$, where $\nu,\ga\ge 1$,
and the constant $c$ does not depend on $Q$. Then for any $\nu+1\le
p\le r<\infty$ and for all $f$ with $(Tf)^*(+\infty)=0$,
\begin{equation}\label{ineq}
\|Tf\|_{L^p(w)}\le
c\|w\|_{(A_p)^{\frac{1}{p-1}}(A_r)^{\frac{1}{\nu}-\frac{1}{p-1}}}\|f\|_{L^p(w)},
\end{equation}
where $c=c(T,p,r,\nu,\ga,n).$
\end{theorem}

If it is known additionally that $T$ is, for example, of weak type
$(1,1)$ (which is the case of any operator from Theorems \ref{mixcz} and \ref{mixlp}),
then $(Tf)^*(+\infty)=0$ for any $f\in L^1$. Hence, we get
first (\ref{ineq}) for $f\in L^1\cap L^p(w)$, and then by the
standard argument it is extended to any $f\in L^p(w)$.

Condition (\ref{cond}) for the maximal Haar shift operator ${\mathbb
S}^{*}f$ was proved in \cite{CMP} (see also \cite{L}) with $\nu=1$
and $\ga$ depending on the complexity. Hence, by the above
discussion in Section 2.1, Theorem \ref{key} implies Theorem
\ref{mixcz}.

Further, in the case $\nu=2$ condition (\ref{cond}) holds for the
dyadic square function $S_d$ with $\ga=1$ (this fact was proved in
\cite{CMP}), and for the intrinsic square function $\widetilde
G_{\a}$ with $\ga=15$ (this was proved in \cite{L2}). From this and
from Theorem \ref{key} we get Theorem \ref{mixlp}.

\begin{proof}[Proof of Theorem \ref{key}]
Combining (\ref{cond}) with Theorem \ref{decom}, we get that for
a.e. $x\in Q_0$,
$$||Tf(x)|^{\nu}-m_{|Tf|^{\nu}}(Q_0)|^{1/\nu}\le c\big(Mf(x)+{\mathcal A}_{3\ga,\nu}f(x)\big),$$
where
$${\mathcal A}_{\ga,\nu}f(x)=\left(\sum_{j,k}\Big(\frac{1}{|\ga Q_j^k|}\int_{\ga
Q_j^k}|f|dx\Big)^{\nu}\chi_{Q_j^k}(x)\right)^{1/\nu}.$$

Therefore, the proof will follow from the corresponding bounds for
$M$ and ${\mathcal A}_{\ga,\nu}$. After that, letting $Q_0$ to anyone
of $2^n$ quadrants, we get that $m_{|Tf|^{\nu}}(Q_0)\to 0$ (since
$(Tf)^*(+\infty)=0$), and Fatou's theorem would complete the proof.

By Buckley's theorem \cite{B},
$\|M\|_{L^p(w)}\le c(p,n)\|w\|_{A_p}^{\frac{1}{p-1}}$, which implies trivially the desired bound for $M$.
Therefore, the proof is reduced to showing that for any $\nu+1\le p\le r<\infty$,
\begin{equation}\label{ag}
\|{\mathcal A}_{\ga,\nu}f\|_{L^p(w)}\le
c\|w\|_{(A_p)^{\frac{1}{p-1}}(A_r)^{\frac{1}{\nu}-\frac{1}{p-1}}}\|f\|_{L^p(w)},
\end{equation}
where $c=c(p,r,\nu,\ga,n).$

In order to handle ${\mathcal A}_{\ga,\nu}f$, following \cite{CMP},
we use the duality. There exists a function $h\ge 0$ with
$\|h\|_{L^{(p/\nu)'}(w)}=1$ such that
$$\|{\mathcal A}_{\ga,\nu}f\|_{L^p(w)}=
\|{\mathcal A}_{\ga,\nu}f\|_{L^{\nu}(hw)}.$$ Further,
\begin{eqnarray}
&&\int_{{\mathbb R}^n}({\mathcal A}_{\ga,\nu}f)^{\nu}hw=\sum_{j,k}
\Big(\frac{1}{|\ga Q_j^k|}\int_{\ga
Q_j^k}|f|dx\Big)^{\nu}\int_{Q_j^k}hw\label{equ}\\
&&=\sum_{j,k}\Big(\frac{\si(3\ga Q_j^k)}{|\ga
Q_j^k|}\Big)^{\nu}\Big(\frac{w(Q_j^k)}{|Q_j^k|}\Big)^{\frac{\nu}{p-1}}
\Big(\frac{1}{\si(3\ga Q_j^k)}\int_{\ga
Q_j^k}|f|dx\Big)^{\nu}\nonumber\\
&&\times\Big(\frac{1}{w(Q_j^k)}\int_{Q_j^k}hw\Big)|Q_j^k|^{\frac{\nu}{p-1}}w(Q_j^k)
^{1-\frac{\nu}{p-1}}.\nonumber
\end{eqnarray}

It is well-known that, by H\"older's inequality, $1\le A_r(w;E)$ for any measurable set $E$ with $|E|>0$.
From this, for any $E\subset Q$ with $|E|\ge \xi|Q|$,
\begin{eqnarray*}
w(Q)&\le&
w(Q)\frac{1}{|E|^r}\Big(\int_Ew^{-\frac{1}{r-1}}\Big)^{r-1}w(E)\\
&\le&(|Q|/|E|)^rA_r(w;Q)w(E)\le (1/\xi)^{r}A_r(w;Q)w(E).
\end{eqnarray*}
Therefore, $w(Q_j^k)\le 2^rA_r(w;Q_j^k)w(E_j^k)$ (the sets $E_j^k$ are defined after Theorem \ref{decom}).
Combining this with (\ref{equ}), we get
\begin{eqnarray*}
&&\int_{{\mathbb R}^n}({\mathcal A}_{\ga,\nu}f)^{\nu}hw\\
&&\le 2^{r(1-\frac{\nu}{p-1})}\sum_{j,k}\Big(\frac{\si(3\ga Q_j^k)}{|\ga
Q_j^k|}\Big)^{\nu}\Big(\frac{w(Q_j^k)}{|Q_j^k|}\Big)^{\frac{\nu}{p-1}}
A_r(w;Q_j^k)^{1-\frac{\nu}{p-1}}\\
&&\times \Big(\frac{1}{\si(3\ga Q_j^k)}\int_{\ga
Q_j^k}|f|dx\Big)^{\nu}\Big(\frac{1}{w(Q_j^k)}\int_{Q_j^k}hw\Big)
|Q_j^k|^{\frac{\nu}{p-1}}w(E_j^k) ^{1-\frac{\nu}{p-1}}.
\end{eqnarray*}
Since
\begin{eqnarray*}
&&\Big(\frac{\si(3\ga Q_j^k)}{|\ga
Q_j^k|}\Big)^{\nu}\Big(\frac{w(Q_j^k)}{|Q_j^k|}\Big)^{\frac{\nu}{p-1}}
A_r(w;Q_j^k)^{1-\frac{\nu}{p-1}}\\
&&\le cA_p(w;3\ga Q_j^k)^{\frac{\nu}{p-1}}A_r(w;3\ga
Q_j^k)^{1-\frac{\nu}{p-1}}\le
c\|w\|_{(A_p)^{\frac{\nu}{p-1}}(A_r)^{1-\frac{\nu}{p-1}}},
\end{eqnarray*}
we obtain
\begin{eqnarray*}
&&\int_{{\mathbb R}^n}({\mathcal A}_{\ga,\nu}f)^{\nu}hw \le c
\|w\|_{(A_p)^{\frac{\nu}{p-1}}(A_r)^{1-\frac{\nu}{p-1}}}\\
&&\times \sum_{j,k}\Big(\frac{1}{\si(3\ga Q_j^k)}\int_{\ga
Q_j^k}|f|dx\Big)^{\nu}\Big(\frac{1}{w(Q_j^k)}\int_{Q_j^k}hw\Big)
|Q_j^k|^{\frac{\nu}{p-1}}w(E_j^k) ^{1-\frac{\nu}{p-1}}.
\end{eqnarray*}
By H\"older's inequality,
\begin{eqnarray*}
&&\sum_{j,k}\Big(\frac{1}{\si(3\ga Q_j^k)}\int_{\ga
Q_j^k}|f|dx\Big)^{\nu}\Big(\frac{1}{w(Q_j^k)}\int_{Q_j^k}hw\Big)
|Q_j^k|^{\frac{\nu}{p-1}}w(E_j^k) ^{1-\frac{\nu}{p-1}}\\
&&\le \left(\sum_{j,k}\Big(\frac{1}{\si(3\ga Q_j^k)}\int_{\ga
Q_j^k}|f|dx\Big)^p|Q_j^k|^{\frac{p}{p-1}}w(E_j^k)^{-\frac{1}{p-1}}\right)^{\nu/p}\\
&&\times
\left(\sum_{j,k}\Big(\frac{1}{w(Q_j^k)}\int_{Q_j^k}hw\Big)^{(p/\nu)'}w(E_j^k)\right)^{1-\nu/p}.
\end{eqnarray*}

Let $M_w^c$ and $M_w^d$ be the weighted centered and dyadic maximal
operator, respectively. We will use below the well known fact that
these operators are bounded on $L^p(w),p>1,$ with the corresponding
bounds independent of $w$.

Since $1\le A_p(w;E_j^k)$, we get from this that
$$|Q_j^k|^{\frac{p}{p-1}}w(E_j^k)^{-\frac{1}{p-1}}\le
2^{\frac{p}{p-1}}|E_j^k|^{\frac{p}{p-1}}w(E_j^k)^{-\frac{1}{p-1}}\le
2^{\frac{p}{p-1}}\si(E_j^k),$$ and therefore,
\begin{eqnarray*}
&&\sum_{j,k}\Big(\frac{1}{\si(3\ga Q_j^k)}\int_{\ga
Q_j^k}|f|dx\Big)^p|Q_j^k|^{\frac{p}{p-1}}w(E_j^k)^{-\frac{1}{p-1}}\\
&&\le c\sum_{j,k}\int_{E_j^k}M_{\si}^c(f/\si)^p\si dx\le
c\int_{{\mathbb R}^n}M_{\si}^c(f/\si)^p\si dx\le c\int_{{\mathbb
R}^n}|f|^pw dx.
\end{eqnarray*}
Similarly,
\begin{eqnarray*}
&&\sum_{j,k}\Big(\frac{1}{w(Q_j^k)}\int_{Q_j^k}hw\Big)^{(p/\nu)'}w(E_j^k)\le
\sum_{j,k}\int_{E_j^k}(M^{d}_wh)^{(p/\nu)'}wdx\\
&&\le \int_{{\mathbb R}^n}(M^{d}_wh)^{(p/\nu)'}wdx\le
c\int_{{\mathbb R}^n}h^{(p/\nu)'}wdx=c.
\end{eqnarray*}

Combining the previous estimates yields
$$\int_{{\mathbb R}^n}({\mathcal A}_{\ga,\nu}f)^{\nu}hw \le c
\|w\|_{(A_p)^{\frac{\nu}{p-1}}(A_r)^{1-\frac{\nu}{p-1}}}\Big(\int_{{\mathbb
R}^n}|f|^pw dx\Big)^{\nu/p},$$ which implies (\ref{ag}), and
therefore, the proof is complete.
\end{proof}

\begin{remark}\label{vvm}
Theorem \ref{key} can be also used to get a new bound for the
vector-valued maximal operator $\overline M_q$ defined for
$f=\{f_i\}$, and $q,1<q<1,$ by
$$
\overline M_qf(x)=\left(\sum_{i=1}^{\infty}Mf_i(x)^q\right)^{1/q}.
$$
It was proved in \cite{CMP} that for any $1<p,q<\infty$,
\begin{equation}\label{vv}
\|\overline M_qf\|_{L^p(w)}\le
c\|w\|_{A_p}^{\max(\frac{1}{q},\frac{1}{p-1})}\Big(\int_{{\mathbb
R}^n}\|f(x)\|_{\ell^q}^pwdx\Big)^{1/p}.
\end{equation}
The proof of this inequality is based on the following variant of
(\ref{cond}) for the vector-valued dyadic maximal operator
$\overline M_q^d$ and any dyadic $Q$:
$$
\o_{\la}((\overline M_q^df)^q;Q)\le
c\Big(\frac{1}{|Q|}\int_{Q}\|f(x)\|_{\ell^q}dx\Big)^q.
$$
Therefore, using the same argument as above, we obtain an
improvement of (\ref{vv}) for $q+1<p\le r<\infty$:
$$
\|\overline M_qf\|_{L^p(w)}\le
c\|w\|_{(A_p)^{\frac{1}{p-1}}(A_r)^{\frac{1}{q}-\frac{1}{p-1}}}
\Big(\int_{{\mathbb R}^n}\|f(x)\|_{\ell^q}^pwdx\Big)^{1/p}.
$$
\end{remark}

\section{Examples}
\subsection{The sharpness of the exponent $\frac{1}{p-1}$}
First we note that the exponent $\frac{1}{p-1}$ in Theorem
\ref{mixcz} is sharp in the sense that
$\|w\|_{(A_p)^{\frac{1}{p-1}}(A_r)^{1-\frac{1}{p-1}}}$ cannot be
replaced by $\|w\|_{(A_p)^{\a}(A_r)^{1-\a}}$ for $\a<\frac{1}{p-1}$.
Indeed, it suffices to consider the same example as in~\cite{B}. Let
$T=H$ is the Hilbert transform. Let $w(x)=|x|^{(p-1)(1-\d)}$ and
$f=|x|^{-1+\d}\chi_{[0,1]}$. Then on one hand we have that
$\|H\|_{L^p(w)}\ge c\d^{-1}$, and on the other hand, if $r>p$, then
$\|w\|_{(A_p)^{\a}(A_r)^{1-\a}}\le c\d^{-\a(p-1)}$. Therefore,
$\a\ge \frac{1}{p-1}$.

The same observation applies to Theorem \ref{mixlp}. For instance,
in the case of the dyadic square function, exactly the same example
as above (see \cite{DGPP}) shows that
$\|w\|_{(A_p)^{\frac{1}{p-1}}(A_r)^{\frac{1}{2}-\frac{1}{p-1}}}$
cannot be replaced by $\|w\|_{(A_p)^{\a}(A_r)^{\frac{1}{2}-\a}}$ for
$\a<\frac{1}{p-1}$.

\subsection{A comparison with M. Lacey's bound} Let $p>2$. We
show that the right-hand sides in (\ref{Lac}) and in Theorem
\ref{mixcz} incomparable.

Let $w=t\chi_{[0,1]}+\chi_{{\mathbb R}\setminus [0,1]}$. It is easy
to see that
$$\|w\|_{A_p}\sim
\|w\|_{(A_p)^{\frac{1}{p-1}}(A_r)^{1-\frac{1}{p-1}}}\sim t.$$
Further, it was shown in \cite{HP} that for any measurable set $E$,
\begin{equation}\label{thp}
\|t\chi_{E}+\chi_{{\mathbb R}\setminus E}\|_{A_{\infty}}'\le 4\log
t\quad(t\ge 3).
\end{equation}
Hence, $\|w\|_{A_{\infty}}'\le 4\log t$ and
$$\|\si\|_{A_{\infty}}'=\|t^{\frac{1}{p-1}}\si\|_{A_{\infty}}'
\le \frac{4}{p-1}\log t\quad(t\ge 3^{p-1}).$$ Therefore,
$$
\|w\|_{A_p}^{1/p}\max\Big(\|w\|_{A_{\infty}}')^{1/p'}
,(\|\si\|_{A_{\infty}}')^{1/p}\Big)\le ct^{1/p}(\log t)^{1/p'},
$$
which shows that the right-hand side in (\ref{Lac}) can be arbitrary
smaller than the one in Theorem \ref{mixcz}.

On the other hand, for $N$ big enough let
$$
w(x) = \begin{cases} |x|^{(p-1)(1-\d)}, & x\in[-1,1] \\
|x-N|^{\d-1}, & x\in[N-1,N+1]\\
1,& \mbox{otherwise.}\end{cases}
$$
Then we have $\|w\|_{A_p}\ge c\d^{-(p-1)}$ (take $I=[0,1]$). Also,
$\|w\|_{A_{\infty}}'\ge c\d^{-1}$ (take $I=[N,N+1]$). Therefore,
$$
\|w\|_{A_p}^{1/p}\max\Big(\|w\|_{A_{\infty}}')^{1/p'}
,(\|\si\|_{A_{\infty}}')^{1/p}\Big)\ge c\d^{-2/p'}.
$$
But for $N$ big enough the supremum defining
$\|w\|_{(A_p)^{\frac{1}{p-1}}(A_r)^{1-\frac{1}{p-1}}}$ can attain on
small intervals containing either $0$ or $N$. If $r>p$, then for any
such interval
$$A_p(w;I)^{\frac{1}{p-1}}A_r(w;I)^{1-\frac{1}{p-1}}\le c/\d,$$
and hence $\|w\|_{(A_p)^{\frac{1}{p-1}}(A_r)^{1-\frac{1}{p-1}}}\le
c/\d.$ This shows that the right-hand side in Theorem \ref{mixcz}
can be arbitrary smaller than the one in (\ref{Lac}).

\end{document}